\documentclass[12pt,twoside]{article}

\usepackage{headerfo}
\usepackage{amsfonts}
\usepackage{amsmath}

\newcounter{theorem}
\newcommand{\newsection}[1]{{\setcounter{theorem}{0} \section{#1}}}
\newtheorem{Theorem}{Theorem}[section]

\newtheorem{Proposition}[Theorem]{Proposition}
\newtheorem{Lemma}[Theorem]{Lemma}

\def\eop{{ \vrule height7pt width7pt depth0pt}\par\bigskip}

\newcommand{\R}{\mathbb R}
\newcommand{\C}{\mathbb R}

\newif\ifpdf
\ifx\pdfoutput\undefined
  \pdffalse
\else
  \pdfoutput=1
  \pdftrue
\fi

\usepackage{latexsym}

\chardef\aa=64

\begin{document}

\renewcommand{\author} {L.~Bos \\  
Department of Computer Science\\
University of Verona \\
Verona, Italy\\}

\renewcommand{\title}{A Characterization of Optimal Prediction Measures via $\ell_1$ Minimization}
\newcommand{\stitle}{$\ell_1$ Characterization of Optimal Prediction Measures}

\renewcommand{\date}{\today}
\flushbottom
\setcounter{page}{1}
\pageheaderlinetrue
\oddpageheader{}{\stitle}{\thepage}
\evenpageheader{\thepage}{\stitle}{}
\thispagestyle{empty}
\vskip1cm
\begin{center}
\LARGE{\bf \title}
\\[0.7cm]
\large{\author}
\end{center}
\vspace{0.3cm}
\begin{center}
\date
\end{center}
\begin{abstract}
Suppose that $K\subset\C$ is compact and that $z_0\in\C\backslash K$ is an external point.  An optimal prediction measure for regression by polynomials of degree at most $n,$ is one for which the variance of the prediction at $z_0$ is as small as possible.  Hoel and Levine (\cite{HL}) have considered the case of $K=[-1,1]$ and $z_0=x_0\in \R\backslash [-1,1],$ where they show that the support of the optimal measure is the $n+1$ extremme points of the Chebyshev polynomial $T_n(x)$ and characterizing the optimal weights  in terms of absolute values of fundamental interpolating Lagrange polynomials.  More recently,  \cite{BLO}  has given the equivalence of the optimal prediction problem with that of finding polynomials of extremal growth. They also study in detail the case of $K=[-1,1]$ and $z_0=ia\in i\R,$ purely imaginary. 
In this work we generalize the Hoel-Levine formula to the general case when the support of the optimal measure is a finite set and give a formula for the optimal weights in terms of a $\ell_1$ minimization problem.
\end{abstract}

\setcounter{section}{0}


\newsection{Introduction} 

\medskip

Optimal Experimental Design has a rich history within Statistics. The interested reader may consult the classical book of
Karlin and Studden \cite{KS} (especially Chapter X) or the more recent monograph of Dette and Studden \cite{DS}.   A brief description of the statistical motivation of the subject of this work is also available in \cite{BLO}.

We state a generalized version of the problem thusly.  Suppose that $A,B\subset \C^d$ are two compact sets.  Let
\[{\cal M}(A):=\{\mu\,:\, \mu \,\,\hbox{is a probability measure on}\,\, A\}\] 
and consider,  for $\{q_1,\cdots,q_N\}\subset \R_n[z],$ a $\mu$-orthonormal basis for $\C_n[z],$
\[ K_n^\mu(z,z):=\sum_{k=1}^N |q_k(z)|^2, \] 
i.e. , the (reciprocal of) the Christoffel function (also known as the Bergman kernel) for $\C_n[z].$

Here
\[N:={\rm dim}(\C_n[z].\]

The probability measure supported on $A,$
\[\mu_n={\rm argmin}_{\mu \in {\cal M}(A)} \max _{z\in B} K_n(z,z)\]
is said to be G-optimal of degree $n$ for the pair $A,B.$

The case of $B=A$ is the classic case of optimal design,  in which case the celebrated equivalence theorem of Kiefer and Wolfowitz \cite{KW}  informs us that G-optimality is equivalent to what is called D-optimality.  Moreover,  in this case,  it was shown in \cite{BBLW} that
the weak-* limit of the optimal $\mu_n$ exists and is equal to the equilibrium measure of complex Pluripotential Theory.

The case of $A=[-1,1]$ and $B=\{x_0\}\subset \R\backslash A,$ a single point,  was considered by Hoel and Levine \cite{HL} and generalized to what we called an optimal prediction measure in \cite{BLO} for 
$A\subset\C^d$ and $B=\{z_0\}\in \C^d\backslash A.$ For simplicity's sake we will often shorten Optimal Prediction Measure to Optimal Measure.

\section{The Case of $A=[-1,1]$ and $B=\{z_0\}\in\R\backslash[-1,1]$}

In \cite{HL} Hoel and Levine show that in the univariate case, for $A=[-1,1],$ and {\it any} $z_0\in \R\backslash A,$  a {\it real} external point, the optimal prediction measure is a discrete measure supported at the $n+1$ extremal points $x_k=\cos(k\pi/n),$ $0\le k\le n,$ of $T_n(x)$ the classical Chebyshev polynomial of the first kind with weights given (also for $z_0\in\C\backslash[-1,1]$) by 

\begin{Lemma} \label{HLWeights} (Hoel-Levine \cite{HL}) Suppose that $-1=x_0<x_1<\cdots<x_n=+1$ are given and that $z_0\in\C\backslash[-1,1].$ Then among all discrete probability measures supported at these points, the measure $\mu=\sum_{i=0}^n w_i\delta_{x_i}$ with
\begin{equation}\label{OptWeights}
w_i:=\frac{|\ell_i(z_0)|}{\sum_{i=0}^n |\ell_i(z_0)|},\,\,0\le i\le n
\end{equation}
with $\ell_i(z)$ the $i$th fundamental Lagrange interpolating polynomial for these points, minimizes
$K_n^\mu(z_0,z_0).$
\end{Lemma}


We remark that in this case it turns out that
\begin{equation}\label{HLthm}
K_n^{\mu_0}(z_0,z_0)=T_n^2(z_0).
\end{equation}

There is also a matrix-vector interpretation of the Hoel-Levine formula that will be useful for our generalization.  

Choose a basis $\{ P_0,\cdots,P_{n}\}$ for the univariate polynomials of degree at most $n$ and let $V\in \R^{(n+1)\times(n+1)}$ be the so-called Vandermonde matrix with components
\[V_{ij}=P_j(x_i),\quad 0\le i,j\le n.\] 
Further, let ${\mathbf p}\in \R^{n+1}$ be the vector of polynomial values 
\[{p}_i=P_i(z_0),\quad 0\le i\le n.\]
The linear system $V^t{\mathbf c}={\mathbf p}$ consists of rows
\[\sum_{j=0}^n (V^t)_{ij} {c}_j ={p}_i=P_i(z_0),\,\, 0\le i\le n,\]
i.e.,
\[\sum_{j=0}^n P_i(x_j){c}_j=P_i(z_0),\,\,0\le i\le n.\]
It follows that
\[ {c}_i=\ell_i(z_0)\]
where $\ell_i(x)$ is the $i$th fundamental Lagrange polynomial for the points $x_0,\cdots,x_n.$

In other words the Hoel-Levine weights are given by
\[w_i=\frac{|{c}_i|}{\sum_{j=0}^n|c_j|}, \quad V^t{\mathbf c}={\mathbf p}.\]

\section{A Generalized Hoel-Levine Formula}

Suppose that $A\subset \C^d$ is compact and that $B\subset \C^d$ consists of a single point $z_0\in \C^d\backslash A.$

We first consider the interpolation case, i.e., when
\[ A=\{ x_1,\cdots,x_N\}\subset \C^d\]
is a discrete set of $N={\rm dim}(\C_n[z])$ distinct points.

\begin{Lemma} \label{HLWeightsd}  Suppose that the set of $N$ distinct points  $A,$ as described above,  are unisolvent for polynomial interpolation, i.e., for every set of values $y_1,\cdots,y_N\in \C$  there exists a unique polynomial $p\in \C_n[z]$ such that
\[p(x_i)=y_i,\quad 1\le i\le N.\]
Suppose further that
 $z_0\in \C^d\backslash A.$ 
 
Then among all discrete probability measures supported at these points, the measure $\mu=\sum_{i=1}^N w_i\delta_{x_i}$ with
\begin{equation}\label{OptWeights2}
w_i:=\frac{|\ell_i(z_0)|}{\sum_{i=1}^N |\ell_i(z_0)|},\,\,1\le i\le N
\end{equation}
with $\ell_i(x)$ the $i$th fundamental Lagrange interpolating polynomial for these points, minimizes
$K_n^\mu(z_0,z_0).$
\end{Lemma}

\noindent {\bf Proof.} We remark that this includes, in particular,  the Hoel-Levine case, Lemma \ref{HLWeights}.

We first note that for such a discrete measure, $\{\ell_i(x)/\sqrt{w_i}\}_{1\le i \le N}$ form an orthonormal basis. Hence
\begin{equation}\label{GenWeights2}
K_n^\mu(z_0,z_0)=\sum_{i=1}^N \frac{|\ell_i(z_0)|^2}{w_i}.
\end{equation}
In the case of the weights chosen according to (\ref{OptWeights2}) we easily obtain
\begin{equation}\label{OptKn}
K_n^{\mu_0}(z_0,z_0)=\left(\sum_{i=0}^n |\ell_i(z_0)|\right)^2.
\end{equation}
We claim that for any choice of weights $K_n$ given by (\ref{GenWeights2}) is at least as large as that given by (\ref{OptKn}). To see this, just note that by the Cauchy-Schwarz inequality,
\begin{eqnarray*}
\left(\sum_{i=1}^N |\ell_i(z_0)|\right)^2&=&\left(\sum_{i=1}^N \frac{|\ell_i(z_0)|}{\sqrt{w_i}}\cdot
\sqrt{w_i}\right)^2\cr
&\le& \left(\sum_{i=1}^N \frac{|\ell_i(z_0)|^2}{w_i}\right)\cdot\left(\sum_{i=1}^N w_i\right)\cr
&=&\sum_{i=1}^N \frac{|\ell_i(z_0)|^2}{w_i}.
\end{eqnarray*}
\noindent \eop

\medskip

Unlike the univariate case, in the multivariate case degeneracies may occur. Here is a simple example. 

Consider the three points (vertices of the standard triangle) in
$\R^2,$
\[A=\{(0,0),(0,1),(1,0)\}\]
and an exterior point,  $(x_0,y_0)$ such that
\[x_0,y_0\ge 0\,\,{\rm and}\,\, x_0+y_0\ge1.\]
The Lagrange polynomials for the points in the order given are
\begin{eqnarray*}
\ell_1(x,y)&=&1-x-y,\cr
\ell_2(x,y)&=&y,\cr
\ell_3(x,y)&=&x.
\end{eqnarray*}
Hence the Lebesgue function
\begin{eqnarray*}
\Lambda&:=&\sum_{i=1}^3 |\ell_i(x_0,y_0)|\cr
&=&|1-x_0-y_0|+|y_0|+|x_0|\cr
&=&(x_0+y_0)-1+y_0+x_0\cr
&=&2(x_0+y_0)-1
\end{eqnarray*}
and the optimal weights, as given by Lemma \ref{HLWeightsd} are
\begin{eqnarray*}
w_1&=&\frac{|1-x_0-y_0|}{\Lambda}=\frac{x_0+y_0-1}{\Lambda},\cr
w_2&=&\frac{|y_0|}{\Lambda}=\frac{y_0}{\Lambda},\cr
w_3&=&\frac{|x_0|}{\Lambda}=\frac{x_0}{\Lambda}.
\end{eqnarray*}
Notice that if the exterior point is on one of the extended edges of the triangle, $x+y=1,$ $y=0,$ or $x=0,$ one of the weights is zero and the measure is degenerate, i.e.,  there is a polynomial $p$ of degree one, not identically zero, for which the integral of $p^2$ is zero.

\bigskip
Now, let $X\subset A$ be a {\it finite} set that {\it contains} the support of an optimal measure $\mu_n,$ which (for the time being) we will assume to be non-degenerate. Note that if $A$ is itself finite then $X=A$ is allowed.

\begin{Theorem} \label{FiniteCase} Suppose that $N:={\rm dim}(\C_n(A))$ and that
$\{P_1,P_2,\cdots,P_N\}$ is a basis for $\C_n(A).$ Further suppose  that $X=\{x_1,\cdots x_M\}$ with $M\ge N.$ Let $V\in C^{M\times N}$ be the rectangular Vandermonde matrix with components
\begin{equation}\label{vdm} 
V_{ij}=P_j(x_i),\quad 1\le i\le M, \,\, 1\le j\le N.
\end{equation}
We assume that $V$ is of full rank $N.$

Further let ${\mathbf p}\in \C^N$ given by
\[p_i:=P_i(z_0).\]
Suppose that 
\[\mu_n:= \sum_{i=1}^M w_i \delta_{x_i}\]
is an optimal measure for degree $n.$ Then there exists a vector ${\mathbf c}\in\R^M$ such that the weight vector
\[{\mathbf w}=\frac{|{\mathbf c}|}{\|{\mathbf c}\|_1}\]
and that
\[{\mathbf c}:={\rm argmin} \{\|{\bf c}\|_1\,:\, V^t{\mathbf c}={\mathbf p}\}.\]

Conversely,  if 
\[{\mathbf c}:={\rm argmin} \{\|{\bf c}\|_1\,:\, V^t{\mathbf c}={\mathbf p}\}\]
and 
\[{\mathbf w}:=\frac{|{\mathbf c}|}{\|{\mathbf c}\|_1},\]
then
\[\mu_n:= \sum_{i=1}^M w_i \delta_{x_i}\]
is an optimal measure for degree $n.$
\end{Theorem}
\noindent {\bf Proof}.  First consider ${\mathbf c}\in \C^M$ such that $V^t{\mathbf c}={\mathbf p}.$ The $i$th row of this system is
\[ \sum_{j=1}^M P_i(x_j)c_j=P_i(z_0).\]
Hence $V^t{\mathbf c}={\mathbf p}$ is equivalent to
\[ \sum_{j=1}^M c_j P(x_j)=P(z_0),\,\, \forall\, P\in \C_n(A).\]
It follows that, for all $P\in \C_n(A),$
\begin{align*}
|P(z_0)|&\le \sum_{j=1}^M |c_j|\, |P(x_j)|\cr
&\le \Bigl( \max_{x\in X} |P(x)|\Bigr) \, \sum_{j=1}^M |c_j|\cr
&= \|P\|_X \, \|{\mathbf c}\|_1.
\end{align*}
Consequently, for every ${\mathbf c}\in \C^M$ such that $V^t{\mathbf c}={\mathbf p},$
\[\|{\mathbf c}\|_1\ge \max_{P\in \C_n(A)} \frac{|P(z_0)|}{\|P\|_X},\]
i.e., the maximal factor of polynomial growth at $z_0$ relative to $X.$ In \cite{BLO} it is shown that the problems of maximal polynomial growth and of minimal variance are equivalent.  Specifically, it is shown that for an optimal prediction measure $\mu_n,$
\[ K_n^{\mu_n}(z_0,z_0)=\left(\max_{P\in \C_n(A)} \frac{|P(z_0)|}{\|P\|_X}\right)^2\]
and hence,  for every ${\mathbf c}\in \C^M$ such that $V^t{\mathbf c}={\mathbf p},$
\begin{equation}\label{eq1}
\|{\mathbf c}\|_1\ge \sqrt{K_n^{\mu_n}(z_0,z_0)}.
\end{equation}

Also,  we will denote a polynomial of {\it extremal growth} by $Q_n(x),$ i.e., $Q_n(x)$ is a polynomial of degree at most $n$ such that $\|Q_n\|_X=1$ and
\[
\max_{P\in \C_n(A)} \frac{|P(z_0)|}{\|P\|_X}=|Q_n(z_0)|.
\]

Now suppose that
\[\mu_n:= \sum_{i=1}^M w_i \delta_{x_i}\]
is an optimal measure for degree $n.$ We first show that there then  exists a vector ${\mathbf c}\in\R^M$ such that the weight vector
\[{\mathbf w}=\frac{|{\mathbf c}|}{\|{\mathbf c}\|_1}\]
and that
\[{\mathbf c}:={\rm argmin} \{\|{\bf c}\|_1\,:\, V^t{\mathbf c}={\mathbf p}\}.\]

To show that our ${\mathbf c}$ has minimal $\ell_1$ norm we will show that it attains inequality (\ref{eq1}). 


\medskip  Now, consider discrete probability measures supported on a subset of $X$ of the form
\[\nu =\sum_{j=1}^M w_j \delta_{x_j}\]
with $w_j\ge0$ and $\sum_{j=1}^M w_j=1.$ For such a measure its  Gram matrix may be written as
\[G=V^tWV\]
where the Vandermonde matrix $V\in \C^{M\times N}$ is defined in (\ref{vdm}) and $W\in \R^{M\times M}$ is the diagonal matrix of the weights, $W_{jj}=w_j,$ $1\le j\le M.$ The kernel is then
\[K({\mathbf w}):=K_n^\nu(z_0,z_0)={\mathbf p}^t G^{-1}{\mathbf p}\]
(where we have made a slight abuse of notation).

We may compute
\begin{align*}
\frac{\partial K}{\partial w_k}&= {\mathbf p}^t \left(-G^{-1}\frac{\partial G}{\partial w_k}G^{-1}\right){\mathbf p}\cr
&=-{\mathbf p}^tG^{-1}(V^t\frac{\partial W}{\partial w_k}V)G^{-1}{\mathbf p}\cr
&={\mathbf p}^tG^{-1}(V^tE_kV)G^{-1}{\mathbf p}
\end{align*}
where  $E_k\in \R^{M\times M}$ is the diagonal matrix with 1 in the $kk$ entry and zero elsewhere.

But it is easy to check that
\[(V^tE_kV)_{ij}=V_{ki}V_{kj},\]
i.e., that
\[V^tE_kV={\mathbf R}_k{\mathbf R}_k^t\]
where ${\mathbf R}_k\in \C^N$ is the $k$th column of $V^t.$

Hence
\[\frac{\partial K}{\partial w_k}=-{\mathbf p}^tG^{-1}{\mathbf R}_k{\mathbf R}_k^tG^{-1}{\mathbf p}.\]
Now, for optimality, either $w_k=0$ or $0<w_k<1$ ($w_k=1$ implies that all the other weights are 0 and hence the measure is supported at a single point, which is not possibly optimal). Thus applying Lagrange multipliers to the constraint $\sum_{j=1}^Mw_j=1$ we must have
\begin{align*}
\frac{\partial K}{\partial w_k}&=-\lambda, \quad w_k\neq0\cr
\iff {\mathbf p}^tG^{-1}{\mathbf R}_k{\mathbf R}_k^tG^{-1}{\mathbf p}&=\lambda,\quad w_k\neq 0\cr
\iff ({\mathbf R}_k^tG^{-1}{\mathbf p})^t({\mathbf R}_k^tG^{-1}{\mathbf p})&=\lambda,\quad w_k\neq 0\cr
\iff ({\mathbf R}_k^tG^{-1}{\mathbf p})^2&=\lambda, \quad w_k\neq0\cr
\iff {\mathbf R}_k^tG^{-1}{\mathbf p}&=\pm \sqrt{\lambda}, \quad w_k\neq0.
\end{align*}
Now,  $K({\mathbf w})$ is homogeneous of order $-1,$ i.e.,
\begin{align*}
K(t{\mathbf w})&={\mathbf p}^t G^{-1}(t{\mathbf w}){\mathbf p}\cr
&= {\mathbf p}^t(V^t(tW)V)^{-1}{\mathbf p}\cr
&=\frac{1}{t} {\mathbf p}^t(V^tWV)^{-1}{\mathbf p}\cr
&=\frac{1}{t} K({\mathbf w}).
\end{align*}
Consequently,  by differentiating with respect to $t$ and setting $t=1$ we obtain
\[\sum_{k=1}^M w_k\frac{\partial K}{\partial w_k}=-K({\mathbf w}).\]
Hence $\partial K/\partial w_k =-\lambda,$ $w_k\neq 0,$ implies that 
\[-\lambda \sum_{k=1}^M w_k=-K({\mathbf w}),\]
i.e., 
\[\lambda = K({\mathbf w}).\]
Therefore the optimality condition is
\begin{equation}\label{optimality}
{\mathbf R}_k^tG^{-1}{\mathbf p}=\pm \sqrt{K({\mathbf w})},\quad w_k\neq 0.
\end{equation}
Vectorizing these equations we obtain
\[VG^{-1}{\mathbf p}=\sqrt{K({\mathbf w})}\,{\mathbf s}\]
where ${\mathbf s}\in \R^M$ is defined as
\[s_k:=({\mathbf R}_k^tG^{-1}{\mathbf p})/\sqrt{K({\mathbf w})},\quad 1\le k\le M.\]
Note that $s_k=\pm 1$ for $w_k\neq0.$

Multiplying on the left by $V^tW$ we obtain
\begin{align*}
(V^tWV)G^{-1}{\mathbf p}&=\sqrt{K({\mathbf w})}\,V^tW{\mathbf s}\cr
&=:V^t{\mathbf c}
\end{align*}
where 
\begin{align*}
{\mathbf c}&:= \sqrt{K({\mathbf w})}\,W{\mathbf s}.
\end{align*}
Note that, as $s_k=\pm 1$ for $w_k\neq 0,$
\[ (W{\mathbf s})_k=\pm w_k,\quad 1\le k\le M.\]
Hence
\[V^t{\mathbf c}=(V^tWV)G^{-1}{\mathbf p}=GG^{-1}{\mathbf p}={\mathbf p}\]
and
\[ \|{\mathbf c}\|_1=\sqrt{K({\mathbf w})}\,\sum_{k=1}^M w_k=
\sqrt{K({\mathbf w})}.\]
Consequently we have proved the first statement of the Theorem.

\medskip

Conversely, suppose now that 
\[{\mathbf c}:={\rm argmin} \{\|{\bf c}\|_1\,:\, V^t{\mathbf c}={\mathbf p}\}\]
and set
\[{\mathbf w}:=\frac{|{\mathbf c}|}{\|{\mathbf c}\|_1}.\]
We will show that
\[\mu_n:= \sum_{i=1}^M w_i \delta_{x_i}\]
is an optimal prediction measure for degree $n.$

It suffices to show that
\[K_n^{\mu_n}(z_0)=|Q_n(z_0)|^2,\]
the polynomial of extremal growth.

Now to see this we will make use of (standard) $\ell_1$ minimization duality.

\begin{Proposition}  
Suppose that $A\in \R^{N\times M}$ and $b\in \R^N.$ Then
\[
\min_{A{\mathbf c}={\mathbf b}} \|{\mathbf c}\|_1=\max_{\|A^t{\mathbf z}\|_\infty\le 1}{\mathbf z}^t{\mathbf b}\]
and, for the optimal vectors,
\[(A^t{\mathbf z})_i={\rm sgn}(c_i),\quad \forall c_i\neq0.\]
\end{Proposition}

In our case we take $A=V^t$ and ${\mathbf b}={\mathbf p}.$ Note that
\[{\mathbf z}^t{\mathbf p}=\sum_{i=1}^N z_iP_i(z_0)\]
is the value of the polynomial, with coefficient vector ${\mathbf z}\in\R^N,$ evaluated at the external point.  
Further, if
\[P(x):=\sum_{i=1}^N z_iP_i(x)\]then
\[V{\mathbf z}=\left[\begin{array}{c}
P(x_1)\cr P(x_2)\cr \cdot\cr\cdot\cr P(x_N)\end{array}\right]\]
and
\[\|V{\mathbf z}\|_\infty=\|P\|_X.\]
Hence
\begin{eqnarray*}
\max_{\|V{\mathbf z}\|_\infty\le 1}{\mathbf z}^t{\mathbf p}
&=& \max_{{\rm deg}(P)\le n,\,\|P\|_X\le 1} |P(z_0)|\cr
&=& \max_{{\rm deg}(P)\le n,\,\|P\|_A\le 1} |P(z_0)|\cr
&=&|Q_n(z_0)|
\end{eqnarray*}
(as we have assumed that $X$ contains the support of an optimal measure).

Consequently,
\[\min_{V^t{\mathbf c}={\mathbf p}} \|{\mathbf c}\|_1=|Q_n(z_0)|.\]

Now to show that
\[K_n^{\mu_n}(z_0)=|Q_n(z_0)|^2\]
we first note that, by the $\ell_1$ duality theorem,  
\[(V{\mathbf z})_i={\rm sgn}(c_i),\quad \forall c_i\neq0.\]
Hence, with ${\mathbf w}:=|{\mathbf c}|/\|{\mathbf c}\|_1$
and
\[W:={\rm diag}({\mathbf w})\in \R^{N\times N}\]
the diagonal matrix with $W_{ii}:=w_i,$ we have
\[{\mathbf c}=\|{\mathbf c}\|_1\times WV{\mathbf z}.\]

Thus
\begin{eqnarray*}
{\mathbf p}&=& |Q_n(z_0)|\,\frac{{\mathbf p}}{\|{\mathbf c}\|_1}\cr\cr
&=& |Q_n(z_0)|\,\frac{V^t{\mathbf c}}{\|{\mathbf c}\|_1}\cr\cr
&=& |Q_n(z_0)|\, V^t(WV{\mathbf z})\cr\cr
&=& |Q_n(z_0)|\, (V^tWV){\mathbf z}\cr\cr
&=& |Q_n(z_0)|\, G{\mathbf z}.
\end{eqnarray*}

Consequently,
\begin{eqnarray*}
K_n^{\mu_n}(z_0)&=&{\mathbf p}^tG^{-1}{\mathbf p}\cr
\cr
&=& |Q_n(z_0)|^2\,{\mathbf z}^tG{\mathbf z}.
\end{eqnarray*}

But
\[{\mathbf z}^tG{\mathbf z}=\int_XP^2(x)\,d\mu_n\]
where
\[P(x)=\sum_{i=1}^Nz_iP_i(x)\]
and, as already noted,
\[V{\mathbf z}=\left[\begin{array}{c}
P(x_1)\cr P(x_2)\cr \cdot\cr\cdot\cr P(x_N)\end{array}\right]\]
so that by the duality theorem,
\[P(x_i)={\rm sgn}(c_i)\]
and hence 
\begin{eqnarray*}
K_n^{\mu_n}(z_0)&=&{\mathbf p}^tG^{-1}{\mathbf p}\cr
\cr
&=& |Q_n(z_0)|^2\,{\mathbf z}^tG{\mathbf z}\cr\cr
&=& |Q_n(z_0)|^2\, \int_X P^2(x)\,d\mu_n(x)\cr\cr
&=& |Q_n(z_0)|^2\, \int_X 1\,d\mu_n(x)\cr\cr
&=& |Q_n(z_0)|^2.
\end{eqnarray*}
The result follows. \eop

\section{In Terms of General Measures}

\begin{Theorem}
Optimal Prediction Measures are given as follows.
Let ${\cal S}(A)$ denote the set of signed measures on $A$ whose support contains the support of an optimal measure, and let
\[\nu = {\rm argmin}\,\,\{|\nu|(A)\,:\,\nu\in{\cal S}(A);\,\,\int_A p(x)d\nu=p(z_0),\,\,\forall p\in \R_n(A)\}.\]
Then
\[\mu_n=\frac{|\nu|}{|\nu|(A)}.\]
\end{Theorem}
\noindent {\bf Proof}. Again, first note that if $\nu\in{\cal S}(A)$ is such that $\int_A p(x)d\nu=p(z_0)$ then
\begin{align*}
|p(z_0)|&\le \int_K |p(x)|d|\nu|\cr
&\le \|p\|_K\times |\nu|(A)
\end{align*}
so that
\begin{equation}\label{lb2}
|\nu|(A)\ge \frac{|p(z_0)|}{\|p\|_A},\quad \forall p\in \R_n(A).
\end{equation}
But by the results of \cite{BLO}, 
for an optimal prediction measure $\mu_n,$
\[ K_n^{\mu_n}(z_0,z_0)=\left(\max_{P\in \R_n(A)} \frac{|P(z_0)|}{\|P\|_A}\right)^2\]
and so 
\begin{equation}\label{lb3}
|\nu|(A)\ge \sqrt{ K_n^{\mu_n}(z_0,z_0)}.
\end{equation}
Just as in the discrete case the result will follow by exhibiting a measure for which (\ref{lb3}) is attained.

Indeed,  suppose that $\mu_n$ is an optimal prediction measure for degree $n$ and $z_0\in \R^d\backslash A.$ By Tchakaloff's Theorem (see e.g.  \cite{Pu97}),  there must exist a measure of finite support, $\mu_n'$ say,  which generates the same Gram matrix and hence is also optimal.  The discrete signed measure with the weights $c_j$ given in the finite case, Theorem \ref{FiniteCase}, using $\mu_n'$ attains (\ref{lb3}). \eop


\bigskip
\begin{thebibliography}{9}

\bibitem{BBLW} \textsc{Bloom,  T.,  Bos,  L.,  Levenberg,  N.  and Waldron,  S.} (2010).  \textit{On the Convergence of Optimal Measures},  Constructive Approx., 32,  59 -- 179.

\bibitem{BLO} \textsc{Bos, L.,  Levenberg,  N.  and Ortega-Cerda,  J.  },  (2020). \textit{Optimal Polynomial Prediction Measures and Extremal Polynomial Growth},   Constructive Approximation,
https://doi.org/10.1007/s00365-020-09522-1


\bibitem{DS} Dette, H. and Studden, W.J., {\bf The Theory of Canonical Moments with Applications in Statistics,
Probability and Analysis}, Wiley Interscience, New York, 1997.

\bibitem{HL} \textsc{Hoel, P.G.  and Levine, A.}, \textit{Optimal spacing and weighting in polynomial prediction},  Ann. Math. Statist. 35 (1964), 1553 -- 1560.

\bibitem{KS} Karlin, S. and Studden, W.J.,  {\bf Tchebycheff Systems: With Applications in Analysis and Statistics}, Wiley Interscience, New York, 1966.


\bibitem{Pu97} M. Putinar, {\it A note on Tchakaloff's theorem}, Proc. Amer. Math. Soc. 125 (1997), 2409--2414.



\end{thebibliography}
\end{document}